\documentclass{amsart}

\usepackage[T1]{fontenc}
\usepackage[utf8]{inputenc}
\usepackage{amsmath,amsthm,amsfonts,amscd,amssymb,eucal,latexsym,mathrsfs}
\usepackage{stmaryrd}
\usepackage{enumerate}
\usepackage{hyperref}
\usepackage[all]{xy}
\usepackage{etoolbox}
\usepackage{amssymb}
\usepackage{amsfonts}
\usepackage{graphicx}
\usepackage{mathtools}

\usepackage{tikz,tikz-cd}
\usetikzlibrary{shapes.geometric}
\usetikzlibrary{arrows}
\usetikzlibrary{calc}

\usetikzlibrary{positioning}
\usepackage{float}
\usepackage{MnSymbol}
\usetikzlibrary{matrix}
\usepackage{tkz-euclide}

\theoremstyle{plain}
\newtheorem{theorem}{Theorem}[section]
\newtheorem*{theorem*}{Theorem}

\newtheorem*{corollary*}{Corollary}

\theoremstyle{definition}

\newtheorem*{definition*}{Definition}

\usetikzlibrary{calc,graphs}
\usepackage{xcolor}
\usetikzlibrary{arrows,decorations.pathmorphing}

\newcommand{\s}{\mathbf{s}}
\renewcommand{\r}{\mathbf{r}}

\newcommand{\ZI}{\mathbb{Z}}

\newcommand{\ts}{\textsection}

\renewcommand{\r}{{\mathbf r}}

\newcommand{\td}{\bigtriangledown}

\newcommand{\B}{\mathcal{B}}

\title{Quasi-periodic tree sets and structure trees}
\author{Mikael Pichot}

\begin{document}

\maketitle
\setcounter{tocdepth}{1}

\begin{abstract}
We expand our discussion of the almost stability theorem of Dicks and Dunwoody in the context of probability measure preserving equivalence relations.
\end{abstract}

\section{Introduction}

This short note is a follow up of our previous note \cite{ast1}. Here we observe that some results by Dunwoody, and by Dicks and Dunwoody, have  analogs in the context of probability measure preserving (pmp) equivalence relations.

\subsection{Tree sets} The following holds.

\begin{theorem}\label{T - Dicks Dunwoody tree set in Bn} 
Let $R$ be a pmp equivalence relation,  $Y$ be a connected $R$-graph, $n\geq 1$ be an integer. There exists a $R$-invariant tree set $E$ which generates  $\B_n(Y)$. 
\end{theorem}

In this theorem, $\B_n(Y)$ denote the set of cuts $s$ in $Y$ whose boundary $\delta s$ contains at most $n$. It is a Boolean ring.

Theorem \ref{T - Dicks Dunwoody tree set in Bn}  result is an direct analog of the following result of Dicks and Dunwoody:

\begin{theorem}[{\cite[Theorem 2.20]{dd}}]\label{T - Dicks Dunwoody tree set in Bn group} 
Let $G$ be a group,  $Y$ be a connected $G$-graph, $n\geq 1$ be an integer. There exists a $G$-invariant tree set $E$ which generates  $\B_n(Y)$. 
\end{theorem}

Theorem \ref{T - Dicks Dunwoody tree set in Bn group} is one of the cornerstone of the Dunwoody's structure tree theory. 
In fact, Dicks and Dunwoody construct an ascending chain $E_1\subset E_2\subset E_3\cdots$ of  tree $G$-sets consisting of thin elements such that 
\begin{enumerate}
\item $E_1$ consists of all $e\mid VY$, where  $e$ is an edge of $Y$ such that $Y\setminus e$ is disconnected;
\item $E_n$  generates $B_n(Y)$ for every $n\geq 1$;
\end{enumerate} 

\noindent (see \cite[p. 61]{dd}). 

A similar statement holds for equivalence relations (see Th.\ \ref{C - Dicks Dunwoody tree set in Bn 2 for equivalence relations}). 

In many contexts, one can construct tree sets directly without relying on the general existence theorem. Thus, the existence of tree sets that are invariant under a group action is a fundamental problem which was famously studied by Dunwoody in connection with his work on accessibility and the Stallings structure theorem \cite[p. 57]{Stallings}. In \cite{Dunwoody-cutting}, Dunwoody constructs tree sets of the form $Gs+Gs^*$, where $G$ is the automorphism group of the graph, and $s$ is a narrow cut (see also   \cite{Moller, Roller, Kron}).   In connection with his accessibility theorem  \cite{Dunwoody-accessible}, he constructs tree sets which are dual to a set of disjoint minimal tracks on a finite simplicial 2-complex  (see also  \cite[Chap.\ IV]{dd}). The general theorem of Dicks and Dunwoody, Theorem \ref{T - Dicks Dunwoody tree set in Bn group}, was first established in their textbook \cite{dd}.

It is not very difficult, but would be quite lengthy, to see that the original proof of this theorem, following \cite[Chap.\ 2]{dd}, extends to equivalence relations. Our original proof of Theorem \ref{T - Dicks Dunwoody tree set in Bn} was written along these lines.  Here we observe that more recent arguments due Dunwoody  \cite{DunwoodyNetworks} make the adaptation of Theorem \ref{T - Dicks Dunwoody tree set in Bn group}  to equivalence relations straightforward (namely, a direct consequence of Duwoody's construction of canonical increasing chains of tree sets). Dunwoody attributes an important idea in his approach to Diekert and Weiss \cite{dd}.  

\subsection{Structure tree}
To every $G$-invariant tree set $E$ is naturally associated a $G$-tree $T(E)$. This tree is called the \emph{structure tree} of $E$. 

The constructions of tree sets, and their associated structure trees, were introduced by Dunwoody in his foundational papers \cite{Dunwoody} (in \ts 2).

 \begin{theorem}[Dunwoody] To every tree set $E$ is naturally associated a tree $T(E)$ with edge set $E$. 
  \end{theorem}

In particular, if $E$ is a $G$-space, then $T(E)$ is a $G$-tree, and if $E$ is an $R$-space, then $T(E)$ (which means applying $T$ fibrewise) is an $R$-tree.

The structure tree $T=T(E)$ of a tree set $E$  is related to the underlying graph $Y$ by  the structure map:
\begin{align*}
\nu\colon VY&\to VT\\
v&\mapsto v \mid E.
\end{align*} 
This map is obtained by viewing a vertex $v\in VY$ as an element in the bidual (i.e., a function of $\B(Y)$), and restrict this element to the subset $E$. We are using here the Dicks--Dunwoody notation (reviewed in \ts\ref{S - Dunwoody} below).

In general, the structure map needs not be injective nor surjective. Its behaviour depends on the geometry of ends in $Y$. If  $\Omega Y$ denote the space of ends of $Y$, then $\nu$ extends naturally to a map:
\[
\nu\colon VY\cup\Omega Y\to VT\cup \Omega T
\]   
which takes vertices to vertices, and may take ends to ends, or to vertices (see \ts\ref{S - Dunwoody} for a description of the extension to ends).

Dunwoody's structure tree construction has immediate corollaries  in the context of equivalence relations.  Thus for instance, if $R$ is an equivalence relation and the tree set $E$ a quasi-periodic in the sense of \cite{qp,ast1} (in the category of countable sets), then $R$ acts naturally on the $R$-tree $T(E)$; in particular, it admits an $R$-tree with a quasi-periodic edge set. Of course, the structure tree  will not be a quasi-periodic tree in general (if the structure tree is quasi-periodic, then $R$ is treeable) without rather strong assumptions on $\nu$; if the structure map $\nu$ is finite-to-1 and $T(E)$ is locally finite, then the structure tree will be a quasi-periodic tree (and therefore $R$ will be treeable). 


We shall only describe the map $\nu$ briefly here. For a more elegant treatment, see for example the paper of Dunwoody \cite{DunwoodyNetworks}. This paper is also (briefly) discussed in \ts\ref{S - Dunwoody} below. Early results on the structure tree construction  include in particular the work by Dunwoody (including \cite{Dunwoody-cutting}), Woess \cite{w} (see in particular in \ts4 therein, but the structure map does not appear explicitly in this reference---see also \ts21 in \cite{wbook}), the work by M\"oller \cite{m1,m2} (see also the survey paper \cite{Moller}), and the work of Thomassen--Woess \cite{tw}.  We will simply quote a few of these standard results here. Suppose for example $E\subset \B_n(Y)$ is obtained from Theorem \ref{T - Dicks Dunwoody tree set in Bn group}. Then, $T$ is a $G$-tree, and by deleting an edge from $T$, one divides $VT$ into two, viewing this edge as a cut of $VY$. The boundary of this cut contains at most $n$ edges. Furthermore, by employing all edges in $T$, together with the operations of union, intersection, and complement, one arrives in this way at every partition of $VY$ that arises by deleting at most $n$ edges (see \cite[Chap.\ 2]{dd}). 
 If in addition $Y$ is locally finite, then under the map $\nu$, every end in  $T^x$ naturally corresponds to a unique end in $Y^x$, and every end in $Y^x$ corresponds to either a vertex or an end of $T^x$. Furthermore, the structure map  $\nu\colon VY\to VT$ is finite-to-1, and $T$ is locally finite if and only if it takes no end of $Y$ to a vertex of $T$  (see  \cite[Lemma 4]{m1}).  Observe that these results also have direct analogs in the case of equivalence relations.  For instance, Moller's theorem for equivalence relations show that if  $Y$ is a locally finite quasi-periodic graph, and $E$ is a quasi-periodic tree set in the Boolean ring of $Y$, then  the structure map $\nu\colon VY\to VT$ is finite-to-1 and $T$ is locally finite if and only if it takes no end of $Y$ to vertex of $T$. This result was generalized to non-locally finite graphs by Dunwoody  in \cite[Theorem 2.1]{Dunwoodypro}, and Dunwoody's results in this paper also have interesting consequences for equivalence relations.

\subsection{Almost stability theorem} The present note is part of an ongoing research project on the almost stability theorem to pmp equivalence relations and its applications \cite{ast1}. The proof of the almost stability theorem in \cite{dd} is divided into two chapters.  The present note only addresses Chap.\ 2 which is the easier part. The author's initial proof of Theorem \ref{T - Dicks Dunwoody tree set in Bn} was a direct modification of the proof given in \cite{dd}.  It is a proof by induction that relies on a maximality argument. This type of argument is standard in the theory of equivalence relations and von Neumann algebras, and  the  proof can be adapted with some (reasonable) amount of effort. The maximality argument in \cite{dd} involves a  specific  choice of thin cuts to enlarge a tree set from step $n$ to $n+1$. A more canonical construction  would simplify the proof for equivalence relations substantially. Dunwoody's paper \cite{Dunwoody} provides such a construction. One of my former students, Jonah Saks, participated in this project. His MSc research topic was to determine which portion of Chap.\ 3 is required to prove the results in Chap. 4.1 (on finite index) of \cite{dd} for equivalence relations. Prior to his MSc, Jonah also did undergraduate research under my supervision, also on  \cite{dd}, and more precisely on the question (which turned out to not be open \cite{Dunwoody}) of eliminating the choice made from step $n$ to $n+1$ in Chap.\ 2 of \cite{dd}.

\subsection{Borel graphs} Important recent works on related subjects take place in the context of Borel graphs \cite{T,Poulain,JS,I,anush}.  In \cite{T}, Tserunyan proves an analog of Stallings theorem for Borel graphs; in \cite{Poulain}, the authors prove that every locally finite whose large scale geometry is tree-like, induces a treeable equivalence relation;  in \cite{JS}, Jardon-Sanchez studies Borel graphs with accessible planar connected components;  in \cite{I}, Ishikura proves an analog for Borel graphs of Dunwoody's 1-cohomology characterization of accessibility; 
in \cite{anush}, the authors propose a direct LLM assisted solution to the finite index treeability question. Some of these papers are written in the formalism of wall spaces, ultrafilters, and CAT(0) cubes complexes, median graphs. The work of Roller \cite{Roller} links Dunwoody's construction to this formalism, and a more precise connection exists in some cases. In \cite{Poulain} for instance, the finiteness requirement on  hyperplanes puts all half-spaces in the Boolean ring of cuts, where Dunwoody's theory could be applied. We will, however, continue to work in the context of Dunwoody's theory in this project (perhaps with the help of an LLM!).

\subsection{On AI usage} No use of AI was made at any stage in the conception or writing of the present note.

\section{Dunwoody's theorem}\label{S - Dunwoody}

  We follow the notation in the textbook \cite{dd} (and our previous note \cite{ast1}).  Thus, if $V$ is a countable set, we  let $(V,\ZI_2)$ denote the power set of $V$, which can be viewed as the set of functions on $V$ taking values in $\ZI_2:=\ZI/2\ZI$. If $a\in (V,\ZI_2)$ we write $a^*=1-a$. If $a,b\in (V,\ZI_2)$ we define $a\square b:=(ab, {a}^*b,a{b}^*,{a}^*{b}^*)$ and write $a\leq b$ if $ab^*=0$. Thus, $0\in a\square b$ if and only if $a\leq b^*$, $b\leq a$, $a\leq b$, or $a^*\leq b$; when these conditions hold, we say that $a$ and $b$ are \emph{nested}. We write $a<b$ if $a\leq b$ but $a\neq b$. 

Let $E$ be a countable set endowed with a map $E\to  (V,\ZI_2)$. If $e\in E$, we write $e\mid V$ for the image of $e$ in $(V,\ZI_2)$, and write $e(v)$ for the value of this function at $v$. The set of all such function is denoted $E\mid V$. Furthermore, we endow the set $V$ with the dual map $V\to (E,\ZI_2)$ taking $v$ to the function $e\mapsto e(v)$, and use the same notation $v\mid E$ and $V\mid E$ for the dual.

If $u,v\in (E,\ZI_2)$, we write $u\td v:=\{e\in E: u(e)\neq v(e)\}$ and say that $u$ and $v$ are \emph{almost equal} if this set is finite.  Almost equality is an equivalence relation on $(E,\ZI_2)$. 

  By \emph{countable graph}, we mean a pair $Y=(V,E)$ of countable sets endowed with two maps $s,r\colon E\to V^2$. A graph morphism means the usual thing. We write $VY=V$ for the vertex set and $EY=E$ for the edge set of $Y$.  A \emph{tree} is a connected acyclic countable graph.

 By (countable) \emph{tree set} we mean, following \cite[Chap.\ II, Def.\ 1.3]{dd},  a (countable) subset $E$ of $(V,\ZI_2)$ such that
 \begin{enumerate}
 \item $E$ is nested;
 \item $E$ does not contain any constant functions;
 \item  $V\mid E$ lies in an almost equality class of $(E,\ZI_2)$.  
 \end{enumerate} 
 Condition (3) means that for all $u,v\in V$, the functions $u\mid E$ and $v\mid E$ agree up to a finite set. Thus, if we view the elements of $E$ as cuts, only finitely many cuts in $E$ will separate $u$ and $v$ (finite interval condition).

Associated with a tree set $E$ is a tree $T=T(E)$ in which $ET=E$ and the attaching maps are defined for every $e,f\in E$ by
\[
\s(e)(f):=\begin{cases} 1 \text{ if } f>e\text{ or  } f>e^*\\ 0\text{ if } f\leq e\text{ or  } f\leq e^*\end{cases} \ \ \text{ and }\ \ \
\r(e)(f):=\begin{cases} 1 \text{ if } f\geq e\text{ or  } f>e^*\\ 0\text{ if } f< e\text{ or  } f\leq e^*\end{cases} 
\]
Thus in this description $\s(e),\r(e)\in (E,\ZI_2)$, where precisely one of the following six  cases occurs: $e=f$, $e=f^*$, $f>e$, $f^*>e$, $f>e^*$, $f^*>e^*$ for every pair $(e,f)\in E$. The vertex set of this tree is $VT:=\s(E)\cup \r(E)$. This definition is from \cite{dd}; see also \cite{Dunwoody}.

\begin{theorem}[{\cite[Chap.\ II]{dd}}]\label{T - realization tree set}
$T$ is a tree.
\end{theorem}

The tree $T=T(E)$ is usually referred to as the \emph{structure tree} associated with $E$ in a situation where $V$ is endowed with a group action and $E$ is invariant. We shall continue to use this terminology in more general contexts, in particular, when $E$ is invariant under the action of an equivalence relation.
Trees sets (viewed as partially ordered sets) and their associated structure trees were   introduced by Dunwoody in \cite[\ts 2]{Dunwoody}. 

Assume for simplicity that $E\cap E^*=\emptyset$. The general case is discussed in \cite{dd}. We extend any $v\in (E,\ZI_2)$ to $E\cup E^*$ by setting $v(e^*):=1-v(e)$ and call the set $E\cup E^*$ an \emph{undirected tree set}. We view $T=T(E\cup E^*)$ as an undirected tree.   A characterization in terms of partially ordered sets is discussed in  \cite{dd} (see Chap.\ II, Cor. 1.10).

Let $E$ be an undirected tree set. An element $v\in (E,\ZI_2)$ is said to \emph{respect the nesting}  if $ef=0$ in implies $v(e)v(f)=0$. Here $v(e^*)=1-v(e)$ for all $e\in E=E^*$. It is shown in the proof of \cite[Theorem 1.5]{dd}  that $VT$ coincides with the set $V(E)$ of elements in the almost equality class of $V\mid E$ that preserve the nesting. In particular, $VT$ contains the image $V\mid E$ of the natural map 
\begin{align*}
V&\to VT\\
v&\mapsto v\mid E
\end{align*}
which restricts an element $v$, viewed as an element in the bidual, to the set $E$.

 Let $Y$ be a connected countable graph. A function $s\in (VY,\ZI_2)$ is said to cut  an edge in $EY$ from $u$ to $v$ if $s(u)\neq s(v)$. We write $\delta s$ for the set of edges that are cut by $s$. We call component of $s$ the connected component $Y\setminus \delta s$. The \emph{Boolean ring} of $Y$ is the set $\B(Y)$ of all $s\in (VY,\ZI_2)$ such that $|\delta s|<\infty$. It is a subring of the Boolean ring $(VY,\ZI_2)$.  A non constant element in $\B(Y)$ is called a cut of $Y$. We call $|\delta s|$ the capacity of a cut $s$.

Let $Y$ be a countable graph, let $\B(Y)$ denote the Boolean ring of $Y$. Dunwoody \cite{DunwoodyNetworks} proves that a canonical increasing chain of tree sets exists in $\B(Y)$. In order to discuss this result, we need to  recall some terminology from \cite{DunwoodyNetworks}.

Let $u,v\in VY\cup \Omega Y$ and  let $s\in \B$. We say that $s$ \emph{separates} $u$ and $v$  if one belongs to $s$ and the other to $s^*$. We say that  $s$ is \emph{thin with respect to $u$ and $v$} if it separates $u$ and $v$ and $|\delta s|$ is minimal among all cuts in $\B$ that separates $u$ and $v$. A cut is \emph{thin} if it is thin with respect to some $u,v$.

This notion of thinness was introduced in \cite{dw}. It is equivalent (see \cite{DunwoodyNetworks}) to the notion used in \cite{dd}, and is useful to construct canonical tree sets. 

A cut $s$ is said to be \emph{tight} both $s$ and $s^*$ are connected \cite{tw}.  Thin cuts are  tight \cite[Prop.\ 2.3]{DunwoodyNetworks} and  the number of tight cuts that cut any given edge of $Y$ is finite \cite[Prop.\ 2.4]{DunwoodyNetworks}.

A \emph{ray} in $Y$ is an infinite sequence of distinct vertices in which two consecutive vertices are adjacent. If $s\in \B(Y)$ and $\gamma$ is a ray, then infinitely many vertices of $\gamma$ belong either  to $s$ or to $s^*$; we say that $\gamma$ belongs to $s$ in the former case and to $s^*$ in the latter. We say that $s$ separate $\gamma$ and $\gamma'$ if one belongs to $s$ and the other to $s^*$. We say that $\gamma$ and $\gamma'$ are equivalent if no cut of $Y$ separate them. An equivalence class of rays is called an end of $Y$. The set of ends is denoted $\Omega Y$.

Let $E\subset \B(Y)$ be an undirected tree set of thin cuts and let $T=T(E)$. The natural map $VY\to VT$ extends to a natural map $\nu\colon VY\cup \Omega Y\to VT\cup \Omega T$. One way to view the extension to ends, following Dunwoody, is as follows. Any edge $e$ in $Y$ maps to a finite path, whose edges are the elements of $E$ cutting $e$. Concatenating these paths, we see that any ray maps to an infinite path in $T$. Then the element $u\in \Omega Y$ maps to a vertex if this path visits some vertex infinitely times, or to an end if it doesn't. Moller's theorem \cite{m1} shows that if the $\nu$ is finite-to-1 on vertices and $T$ is locally finite, then $\nu$ induces a homeomorphism of $\Omega T$ onto $\Omega T$.

If $s\in \B(Y)$ and $E\subset \B(Y)$ is a set of tight cuts, then the number $\mu(s,E)$ of cuts in $E$ which are not nested with $s$ is finite  \cite[Lemma 2.5]{DunwoodyNetworks}.

\begin{theorem}[Dunwoody's theorem]\label{T - Dicks Dunwoody tree set in Bn 2} 
Let $Y$ be a connected graph. Let $C_n$ denote the set of thin cuts of capacity $n$ in $\B(Y)$.  There exists a canonical ascending chain $E_1\subset E_2\subset E_3\cdots$ of  tree sets of $\B(Y)$ such that 
\begin{enumerate}
\item $E_1$ consists of all $e\mid VY$, where  $e$ is an edge of $Y$ such that $Y\setminus e$ is a disconnected subgraph with components in $\B(Y)$;
\item If  $D_n:=\{s\in C_n,\ \mu(s,E_{n-1})=0\}$ then $E_n=E_{n-1}\cup D_n'$, where $D_n'$ consists of all cuts in $s\in D_n$ for which there exist $u,v\in VY\cup \Omega Y$ such that
\begin{enumerate}
\item $s$ is thin with respect to $u$ and $v$;
\item   $\mu(s,D_n)$ is minimal among all values $\mu(t,D_n)$, where $t\in D_n$ is thin with respect to $u$ and $v$;
\end{enumerate} 
when these two conditions hold one says that $s$ is optimally $n$-nested. Thus, $E_n$ is obtained from $E_{n-1}$ by adjoining all optimally $n$-nested cuts.
\item $E_n$ generates $\B_n(Y)$.
\end{enumerate}
In particular, the direct union $E_\omega := \varinjlim E_n$ is  a nested set of thin cuts that generates $\B(Y)$.
\end{theorem}

It following immediately from this result that, if $Y$ is an $G$-graph, $n\geq 1$ is an integer, then there exists a tree $G$-set $E_n$ that generates $\B_n(Y)$. Furthermore, there exists a $G$-tree $T_n$ and a natural $G$-map $VY\to VT_n$ such that $ET_n=E_n$ and $VT_n=V(E_n)$. 

In the context of Borel equivalence relation (or more generally for Borel groupoids),  the following exact analog holds:

\begin{theorem}
\label{C - Dicks Dunwoody tree set in Bn 2 for equivalence relations} 
Let $R$ be an equivalence relation and $Y$ be a connected $R$-graph. Let $C_n$ denote the set of thin cuts with capacity $n$ in $\B(Y)$.  

There exists a canonical ascending chain $E_1\subset E_2\subset E_3\cdots$ of  tree $R$-sets such that 
\begin{enumerate}
\item $E_1^x$ consists of all $e\mid VY^x$, where  $e$ is an edge of $Y^x$ such that $Y^x\setminus e$ is a disconnected and subgraph with components in $\B^x(Y)$;
\item If  $D_n^x:=\{s\in C_n^x,\ \mu(s,E_{n-1}^x)=0\}$ then $E_n^x=E_{n-1}^x\cup D_n'^x$, where $D_n'^x$ consists of all cuts in $s\in D_n^x$ which are optimally nested with respect to some $u,v\in VY^x\cup \Omega Y^x$;
\item $E_n^x$ generates $\B_n^x(Y)$.
\end{enumerate}
In particular, the direct union $E_\omega := \varinjlim E_n$ is then a nested set that generates $\B(Y)$.
\end{theorem}

\begin{proof} Since  Dunwoody's construction is canonical,
it suffices to apply Theorem \ref{T - Dicks Dunwoody tree set in Bn 2} fiberwise. The graph isomorphism $(x,y)\colon Y^y\to Y^x$ provide an (obviously Borel) action on the canonical field of tree sets. For each $x\in X$, the elements of $E_n^x$ satisfy the assumptions of Theorem \ref{T - Dicks Dunwoody tree set in Bn 2} are therefore form a nested set of thin cuts which satisfies (1), (2) and (3). 

It follows by \cite{ast1} that if $EY$ is quasi-periodic, then $E_n$ is quasi-periodic for every $n\geq 1$. In particular, the direct union $E_\omega := \varinjlim E_n$ remains quasi-periodic.
\end{proof}


\begin{thebibliography}{00}


\bibitem{Poulain}
Chen, R., Poulin, A., Tao, R., Tserunyan, A.  "Tree-like graphings, wallings, and median graphings of equivalence relations." Forum of Mathematics, Sigma. Vol. 13. Cambridge University Press, 2025.

\bibitem{anush}
Chen, R.,  Tserunyan, A., Tucker-Drob, R.,  "Finite index extensions of treeable equivalence relations are treeable.", preprint (October 2026).




\bibitem{dicks} Dicks, W., An improved proof of the Almost Stability Theorem (2018).





\bibitem{dd} 
Dicks, W., Dunwoody, M. J., Groups Acting on Graphs.
Cambridge University Press, 1989, 283 pages

\bibitem{dw}
Diekert V., and Weiss, A. Context-free groups and their structure trees. International Journal of Algebra and Computation 23.03 (2013): 611-642.

\bibitem{Dunwoody} 
Dunwoody, M. J., Accessibility and groups of cohomological dimension one, Proc. London Math. Soc. (3), 38(2):193–215, 1979.

\bibitem{Dunwoody-cutting} 
Dunwoody, M. J. Cutting up graphs. Combinatorica 2.1 (1982): 15-23.

\bibitem{Dunwoody-accessible} 
Dunwoody, M. J. The accessibility of finitely presented groups. Invent. math 81.3 (1985): 449-457.

\bibitem{Dunwoodypro}  Dunwoody, M. J., Inaccessible groups and protrees
 Journal of Pure and Applied Algebra 88 (1993) 63-78

\bibitem{DunwoodyNetworks} 
Dunwoody, M. J., Structure trees, networks, and almost invariant sets (2015). 

\bibitem{I} Ishikura H., Decomposition of Borel graphs and cohomology, arXiv, 2026. 

\bibitem{JS}
Jardon-Sanchez H., Applications of tree decompositions and accessibility to treeability of Borel graphs,
preprint, arXiv:2308.13087.

\bibitem{KM} Kechris, A. S., Miller, B. D. Topics in orbit equivalence (No. 1852). Springer Science. (2004)

\bibitem{Kron}
Kr\"on, B., Quasi-isometries between non-locally-finite graphs and structure trees, Abhandlungen aus dem Mathematischen Seminar der Universität Hamburg. Vol. 71. No. 1. Berlin/Heidelberg: Springer-Verlag, 2001.


\bibitem{m1}
Moller, R. G. Ends of graphs. In Math. Proc. Camb. Phil. Soc, Vol. 111, No. 2, pp. 255-266 (1992, March). 

\bibitem{m2}
M\"oller, R. Ends of graphs, II. Math. Proc. Camb. Phil. Soc. (1992), 111, 455. 


\bibitem{Moller}
M\"oller, R. Groups acting on locally finite graphs---a survey of the infinitely ended case. Groups. Vol. 93. 1995.





\bibitem{qp}
Pichot M., Quasi-periodicit\'e et th\'eorie de la mesure, Ph.\ D.\ thesis (2004).

\bibitem{ast1}
Pichot M., Quasi-periodicity and almost equality classes, 2026.



\bibitem{Roller} 
Roller, Martin. ''Poc sets, median algebras and group actions'', arXiv:1607.07747 (1998).

\bibitem{Saks} Saks, J., McGill MSc thesis.

\bibitem{Stallings} 
Stallings, J. R., Group theory and three-dimensional manifolds  (1971).

\bibitem{tw}
Thomassen C., and Woess W., Vertex-transitive graphs and accessibility. Journal of Combinatorial Theory, Series B 58.2 (1993): 248-268.



\bibitem{T} Tserunyan A., A descriptive construction of trees and Stallings’ theorem, Trends in set theory, 191–207,
Contemp. Math., 752, Amer. Math. Soc., Providence, RI, 2020.


\bibitem{w} Woess, W., Boundaries of random walks on graphs and groups with infinitely many ends. Israel J. Math. 68 (1989), 271-301.

\bibitem{wbook}
Woess, W. Random Walks on Infinite Graphs and Groups. Cambridge University Press, 2000. Cambridge Tracts in Mathematics, vol. 138. 
\end{thebibliography}
\end{document}